\documentclass{ifacconf}

\usepackage{graphicx}      
\usepackage{natbib}        
\usepackage{tikz,datatool}
\usepgflibrary{fpu}
\usetikzlibrary{positioning,calc,arrows,shapes,shadows,matrix,backgrounds}

\usepackage{amsmath,amsfonts,mathrsfs}
\usepackage{graphicx,fancybox}
\usepackage{arydshln}
\usepackage{pdfsync,upgreek}
\usepackage{MnSymbol}
\usepackage{dsfont}

\renewcommand{\show}[1]{}

\newcommand{\ot}{\!\otimes\!}

\newcommand{\As}{A}
\newcommand{\Bs}{B}
\newcommand{\Cs}{C}
\newcommand{\Ds}{D}

\newcommand{\Af}{A_\nu}
\newcommand{\Bf}{B_\nu}
\newcommand{\Cf}{C_{k,\nu}}
\newcommand{\Df}{1}

\newcommand{\lap}{\la}
\newcommand{\lad}{\t\la}
\newcommand{\nup}{\nu}
\newcommand{\nud}{{\t\nu}}

\newcommand{\Afp}{A_\nup}
\newcommand{\Bfp}{B_\nup}
\newcommand{\Cfp}{C_{\nup}(\lap)}
\newcommand{\Dfp}{D_{\nup}(\lap)}
\newcommand{\Dfpl}[1]{D_{\nup}^{#1}(\lap)}

\newcommand{\Afd}{A_\nud}
\newcommand{\Bfd}{B_\nud}
\newcommand{\Cfd}{C_{\nud}(\lad)}
\newcommand{\Dfd}{D_{\nud}(\lad)}
\newcommand{\Dfdl}[1]{D_{\nud}^{#1}(\lad)}

\newcommand{\Afk}{A_\psi}
\newcommand{\Bfk}{B_\psi}
\newcommand{\Cfk}{C_\psi}

\newcommand{\Ap}{\Afk}
\newcommand{\Bp}{\Bfk}
\newcommand{\Cp}{\Cfp\ot I_n}
\newcommand{\Dp}{\Dfp\ot I_n}

\newcommand{\Ad}{\Afk}
\newcommand{\Bd}{\Bfk}
\newcommand{\Cd}{\Cfd\ot I_n}
\newcommand{\Dd}{\Dfd\ot I_n}

\newcommand{\Afs}{\c{A}}
\newcommand{\Bfs}{\c{B}}
\newcommand{\Cfs}{\c{C}}
\newcommand{\Dfs}{\c{D}}

\newcommand{\mto}{\rightrightarrows}

\newcommand{\T}{\top}

\newcommand{\si}{\sigma}

\newcommand{\theorem}[1]{\begin{thm}#1\end{thm}}
\newcommand{\definition}[1]{\begin{defi}#1\end{defi}}
\newcommand{\corollary}[1]{\begin{coro}#1\end{coro}}
\newtheorem{coro}[thm]{Corollary}
\newtheorem{defi}[thm]{Definition}

\newcommand{\B}{{\cal B}}
\newcommand{\C}{{\cal C}}

\newcommand{\Dset}{{\Bbb D}}

\newcommand{\z}{{\rm z}}

\newcommand{\eig}{\text{eig}}
\newcommand{\col}{\text{col}}

\newcommand{\te}[1]{\text{\ \ #1\ \ }}

\newcommand{\skipthis}[1]{ }

\renewcommand{\t}[1]{ \tilde{#1} }

\renewcommand{\c}[1]{{\cal #1}}

\newcommand{\eps}{\epsilon}

\newcommand{\R}{ {\mathbb{R}} }

\newcommand{\N}{ {\mathbb{N}} }

\renewcommand{\l}[1]{\label{#1}}
\renewcommand{\r}[1]{(\ref{#1})}

\newcommand{\mun}[1]{\begin{multline*}#1\end{multline*}}
\newcommand{\mul}[1]{\begin{multline}#1\end{multline}}
\newcommand{\equ}[1]{\begin{equation}#1\end{equation}}
\newcommand{\eql}[2]{\begin{equation}\label{#1}#2\end{equation}}

\newcommand{\gan}[1]{\begin{gather*}#1\end{gather*}}

\newcommand{\arr}[2]{\begin{array}{#1}#2\end{array}}

\newcommand{\mat}[2]{\left(\begin{array}{#1}#2\end{array}\right)}
\newcommand{\smat}[2]{{\tiny \left(\begin{array}{#1}#2\end{array}\right)}}

\newcommand{\enu}[1]{\begin{enumerate}#1\end{enumerate}}

\newcommand{\la}{\lambda}
\newcommand{\ga}{\gamma}

\newcommand{\cl}{\prec}
\newcommand{\cg}{\succ}

\newcommand{\cle}{\preceq}

\catcode`\=13
\def{\relax}

\newcommand{\hl}{\\\hline}
\newcommand{\hdl}{\\\hdashline}

\newcommand{\bsub}{\begin{subequations}}
\newcommand{\esub}{\end{subequations}}
\newcommand{\x}{\times}

\newcommand{\ds}{\displaystyle}

\usepackage{enumerate}
\newcommand{\bnumi}{\begin{enumerate}[\upshape (i)]}
\newcommand{\enumi}{\end{enumerate}}

\usetikzlibrary{positioning,calc,arrows,shapes,shadows}
\tikzset{
auto,
sys/.style 2 args={
rectangle,
draw,
drop shadow,
fill=white,
minimum height=#2,
minimum width=#1,
inner sep=\dn},
sum/.style={circle,draw,draw=black,inner sep=0mm,minimum size=2mm},
jun/.style={circle,draw,draw=black,inner sep=0mm,minimum size=0mm},
>={latex},
every path/.style={rounded corners},
}

\def\dn{1ex}
\def\dl{3*\dn}
\def\dh{6*\dn}
\tikzstyle{sy0}=[sys={0*\dn}{0*\dn}]
\tikzstyle{sy1}=[sys={12*\dn}{8*\dn}]
\tikzstyle{sy2}=[sys={8*\dn}{6*\dn}]
\tikzstyle{sy3}=[sys={5*\dn}{5*\dn}]

\newcommand{\fbzf}{
\begin{tikzpicture}[xscale=1,yscale=1]
\tikzstyle{sy4}=[sys={5*\dn}{4*\dn}]
\def\ds{12*\dn}
\def\dl{10*\dn}
\node[sy4] (g) at (0,0)  {$(1-\frac{1}{\z^k})I_d$};
\node[sy4] (d) at ( $(g)-(0,\dh)$ ) {$\partial f$};
\draw[<-] (g)-- node[]{$z$}  ($(g) + (\ds, 0)$);
\draw[->] (g)-- node[swap]{$y$}  ($(g) - (\dl, 0)$);
\draw[->] ($(g) + (.75*\ds, 0)$) |- (d);
\draw[->] (d)-- node[swap,pos=.72]{$w$}  ($(d) - (\dl, 0)$);
\end{tikzpicture}
}

\newcommand{\fbzfa}{
\begin{tikzpicture}[xscale=1,yscale=1]
\tikzstyle{sy4}=[sys={5*\dn}{4*\dn}]
\def\ds{12*\dn}
\def\dl{10*\dn}
\node[sy4] (g) at (0,0)  {$\partial f$};
\node[sy4] (d) at ( $(g)-(0,\dh)$ ) {$(1-\z^k)I_d$};
\draw[<-] (g)-- node[pos=.6]{$z$}  ($(g) + (\ds, 0)$);
\draw[->] (g)-- node[swap,pos=.6]{$w$}  ($(g) - (\ds, 0)$);
\draw[->] ($(g) + (.7*\ds, 0)$) |- (d);
\draw[->] (d)-- node[swap,pos=.4]{$v$}  ($(d) - (\ds, 0)$);
\end{tikzpicture}
}

\begin{document}

\begin{frontmatter}

\title{
Dissipativity, Convexity and Tight O'Shea-Zames-Falb Multipliers for Safety Guarantees
}

\thanks[footnoteinfo]{Funded by Deutsche Forschungsgemeinschaft (DFG, German Research Foundation) under Germany’s Excellence Strategy - EXC 2075 - 390740016. We acknowledge the support by the Stuttgart Center for Simulation Science (SimTech).}

\author[First]{Carsten W. Scherer}

\address[First]{Chair for Mathematical Systems Theory, Department of Mathematics, University of Stuttgart,  Germany.}

\begin{abstract}                
We develop a novel convex parametrization of integral quadratic constraints with a terminal cost for subdifferentials of convex functions, involving general O'Shea-Zames-Falb multipliers. We show the benefit
of our results for the reduction of conservatism of existing techniques, and sketch applications to the analysis of optimization algorithms or the stability analysis of neural network controllers. The development is prepared
by providing a novel link between the convex integrability of a multivariable mapping and dissipativity theory.
\end{abstract}

\begin{keyword}
Dissipativity, Robustness Analysis, Integral quadratic constraints, Absolute Stability, Linear matrix inequalities.
\end{keyword}

\end{frontmatter}

\section{Introduction}\label{Sint}

The field of systems and control is dominated by understanding complex dynamic interconnections. It has been a
highly successful point-of-view to consider such monolithic systems as an interconnection of individual subsystems, and to characterize their global dynamical characteristics by those of the subsystems and properties of their interconnection.
Dissipativity theory as developed by \cite{Wil72a,Wil72b} constitutes a cornerstone in system theory. It not only
translates this conceptual idea into a concrete mathematical framework, but also lays ground for the tailored construction of related computational tools in a modular fashion. For robustness analysis, this theme has been surveyed, e.g., in the recent contributions by \cite{ArcMei16} and \cite{Sch22}, in which one can find many more references to the literature.

The absolute stability analysis of a feedback loop consisting of a linear system and a static nonlinearity, a so-called Lur'e system, has laid ground to this development.
As emphasized in Jan Willems's introduction to the seminal paper by Popov in \cite{Bas01},
the celebrated circle- and Popov-criteria had a substantial impact
on the construction of more flexible and less conservative stability results using
so-called O'Shea-Zames-Falb (OZF) multipliers (\cite{Osh67,WilBro68} and \cite{ZamFal68}).
It is by now well-understood how to generate corresponding computational stability tests
based on the main result in \cite{MegRan97} using integral quadratic constraints (IQCs), which
are strongly inspired by the seminal contributions of \cite{Yak67}. For systems described in
continuous time, computational aspects are exposed in detail in \cite{CheWen96} (see also \cite{VeeSch16}),
while the discrete-time counterparts have been proposed by \cite{CarHea20} and \cite{FetSch17c}.

Most existing stability proofs based on OZF multipliers
are functional analytic in nature. Direct dissipativity-based proofs can be extracted from
\cite{Sei15} and \cite{VeeSch14}.
The latter has lead to the notion of IQCs with a terminal cost, as introduced by \cite{SchVee18}
and further developed by \cite{Sch22} for continuous time systems. It has been argued that this
concept constitutes a seamless link between so-called hard IQCs (which are conservative) and more powerful
soft IQCs on the infinite time-horizon. The papers  by \cite{FetSch17c}, \cite{FetSch17b}, \cite{IanSei19} and \cite{YinSei20} demonstrate the benefit
of such IQCs for reducing conservatism with examples.

Discrete-time absolute stability results based on OZF multipliers have recently drawn considerably attention for the
analysis of optimization algorithms (\cite{LesRec16}), the generalization to the design of extremum controllers (\cite{SchEbe21}), and for the safety verification of neural network controllers (\cite{YinSei20a,PauGra21}). However,
results about IQCs with a terminal cost for discrete-time systems are missing in the literature, which is partly due to technical challenges concerning the correct formulation of the related algebraic Riccati equations.

To fill this gap, this paper develops a novel convex parametrization of IQCs with a nontrivial terminal cost for subdifferentials of convex functions and involving both causal and anti-causal OZF multipliers. We demonstrate the benefit over existing results by a numerical example. It is conjectured that our construction is tight, in the sense
that it constitutes the correct formulation of IQCs in between hard and soft ones; see \cite{FetSch17b} for a discussion of this issue in continuous-time. We only mention various potential applications to optimization algorithms and stability analysis of neural network controllers, the details of which are left for future research.
As a preparation, we reveal a new link between the convex integrability of a multi-valued mapping and dissipativity theory. This is of independent interest on a path towards identifying novel classes of functions for which multiplier-based stability tests can be developed more systematically.

The paper is structured as follows. In Section~\ref{Sdis0}, we relate dissipativity theory to convex integration as
obtained in a celebrated paper by \cite{Roc66}. After recalling strict dissipativity characterizations for
linear systems with quadratic supply rates (Section~\ref{Sdis1}), we formulate a stability result for Lur'e systems
based on IQCs with a terminal cost in discrete-time (Section~\ref{Srob}). Section~\ref{Szf} proposes novel
such IQCs
for subdifferentials of convex functions, and Section~\ref{Salg} reveals how  to bound the amplitude
of an output of a linear system in feedback with any such nonlinearity and in response to initial conditions in the unit ball. A numerical example, a discussion
of the exactness of our main results and of potential applications conclude the paper with
Section~\ref{Sdis}. 

{\bf Notation.} If $X\subset\R^n$, $f:X\mto \R^m$ denotes a possibly multi-valued mapping with domain X. For $m=1$,
the inequality $f(x)\leq 0$ and $\sup_{x\in X} f(x)$ are then interpreted as $y\leq 0$ for all $y\in f(x)$ and
$\sup \{y\ |\ x\in X,\,y\in f(x)\}$, respectively. In $\R^n$, $e$ and $e_k$ are the all-ones and the standard unit vectors for $k\in\{1,\ldots,n\}$. If
$x\in\R^n$, $x\geq 0$ is understood elementwise and $\|x\|$ is the Euclidean norm.
A matrix $A\in\R^{n\times n}$ is Schur if its eigenvalues $\eig(A)$ are located in the unit disk
$\Dset:=\{\la\in\Cset \ |\  |\la|<1\}$. For $A,B\in\R^{n\times n}$, $A\cl B$ means that $A$ and $B$ are symmetric and
$B-A$ is positive definite.
The set of doubly hyperdominant matrices $\c{H}^{n\times n}\subset\R^{n\times n}$ (\cite{WilBro68}) comprises
all $A\in\R^{n\times n}$ whose off-diagonal elements are nonpositive and which satisfy $Ae\geq 0$ and $e^\T A\geq 0$. Note that $\c{H}^{n\times n}$ is spectrahedron, the feasible set of some linear matrix inequality (LMI).
Moreover, $l_{2e}^n$ is the real vector space of all sequences $x:\N_0\to\R^n$, while $l_2^n$ denotes the subspace of $x\in l_{2e}^n$ with $\|x\|_2^2:=\sum_{t=0}^\infty \|x_t\|^2<\infty$.
Finally, we make use of the fact that a linear discrete-time system
$$
x_{t+1}=Ax_t+Bu_t,\ \ y_t=Cx_t+Du_t\te{for}t\in\N
$$
can be lifted as
\eql{lif}{
\mat{cc}{x_{T}\\y^T\\}=\mat{cc}{A^{T}&B^T\\C^T&D^T}\mat{c}{x_0\\u^T}\te{for}T\in\N.
}
This involves stacking the components of the input and output signals $u\in l_{2e}^\bullet$ and $y\in l_{2e}^\bullet$, respectively, as
$$
u^T:=\col\mat{ccccc}{u_0&u_1&\cdots&u_{T-1}},\ \
y^T:=\col\mat{ccccc}{y_0&y_1&\cdots&y_{T-1}},
$$
and using the matrices
$$
A^{T}:=\underbrace{A\,\cdots\,A}_{T\text{\,factors}},\ \ B^T:=\mat{ccccc}{A^{T-1}B&A^{T-2}B&\cdots&B},
$$
$$
C^T :=\mat{c}{C\\CA\\\vdots\\CA^{T-1}},\ \
D^T:=
\mat{ccccc}{
D &0 &\cdots&0\\CB&D&\cdots&0\\\vdots&\ddots&\ddots&\vdots\\CA^{T-2}B&\cdots&CB&D};
$$
for $T=1$ the latter is interpreted as $D^1:=D$.

\section{Dissipativity and Convexity}\l{Sdis0}

For given subsets $X\subset\R^n$, $U\subset\R^m$, $Y\subset\R^k$ and possibly multi-valued
mappings
$$
a:X\times U\mto \R^n,\ \ c:X\times U\mto \R^k,
$$
let us consider the discrete-time dynamical system
\eql{sys}{
\mat{c}{x_{t+1}\\y_t}\in \mat{c}{a(x_t,u_t)\\c(x_t,u_t)}
}
for $t\in\N_0$ with the time axis $\N_0$. We say that the trajectory
$t\mapsto(x_t,u_t,y_t)$ is admissible if
$(x_t,u_t,y_t)\in X\times U\times Y$ holds and if the inclusion \r{sys} is satisfied for all $t\in\N_0$. With the left-shift operator $(\sigma x)(t):=x(t+1)$ for $t\in\N_0$, the system \r{sys} is compactly described as $(\si x,y)\in (a(x,u),c(x,u))$.

\definition{
System \r{sys} is dissipative with respect to the supply rate $S:U\times Y\mto\R$ if there exists a storage function $V:X\to\R$ such that the dissipation inequality (DI)
\eql{di}{
V(x_{t_2})\leq V(x_{t_1})+\sum_{t=t_1}^{t_2-1} S(u_t,y_t).
}
holds for all admissible trajectories and all time instances $t_1,t_2\in\N_0$ with $t_1\leq t_2$. System \r{sys} is cyclo-dissipative
with respect to $S$ in case that
\eql{cdi}{0\leq \sum_{t=t_1}^{t_2-1} S(u_t,y_t)}
holds for all $t_1,t_2\in\N_0$ with $t_1\leq t_2$ and all
admissible trajectories satisfying $x_{t_2}=x_{t_1}$. These properties are referred to as passivity or cyclo-passivity for the particular choice $S(u,y)=u^\T y$ of the supply-rate.
}
Note that, in contrast to the original definition in \cite{Wil72a}, it is not required that storage functions are nonnegative or bounded from below.
Verifying dissipativity requires to come up with a storage function for which the dissipation inequality holds;
we also say that storage functions certify dissipativity.
Instead, cyclo-dissipativity is characterized directly in terms of system trajectories without the need to have a storage function available. Clearly, dissipative systems are cyclo-dissipative.

The converse implication permits  to conclude
the existence of a storage function from cyclo-dissipativity. In other words, a mere trajectory-based input-output
property of the system guarantees the existence of storage functions
which certifies dissipativity.
Under a suitable assumption on the richness of state-trajectories, this converse is established by
explicitly defining two extremal storage functions through the solution of optimal control
problems over admissible system trajectories.

\theorem{\label{Tdis}Suppose there is a ground state $x_*\in X$ such that for every $\xi\in X$ there exists an admissible round trip state trajectory with $x_0=\xi$, $x_{t_1}=x_*$ and $x_{t_2}=\xi$ for some $t_1,t_2\in\N_0$ with $t_1\leq t_2$.
Then \r{sys} is cyclo-dissipative iff it is dissipative as certified by the available storage function
$$
V_*^-(\xi):=\sup_{T\in\N_0,\ x_0=\xi,\ x_{T}=x_*}\left[-\sum_{t=0}^{T-1} S(u_t,y_t)\right]\te{for}\xi\in X,
$$
and the required supply function
$$
V_*^{+}(\xi):=\inf_{T\in\N_0,\ x_0=x_*,\ 
x_{T}=\xi}\left[\sum_{t=0}^{T-1} S(u_t,y_t)\right]\te{for}\xi\in X.
$$
}
This result is proved in complete analogy to the continuous-time counterparts in \cite{HilMoy75} and \cite{Sch21}. The same can be said for the well-known fact that dissipativity is equivalent to the validity of the
local dissipation inequality
\mun{
V(x_+)-V(x)\leq S(u,y)\\\te{for all}(x,u)\in X\times U,\ (x_+,y)\in (a(x,u),c(x,u)).
}

As a first contribution of this paper, let us now establish an interesting relationship of dissipativity with the integration of multi-valued mappings and convexity.
Precisely, we wonder under which  conditions a
mapping $F:\R^n\mto\R^n$ has a convex primitive, i.e., a convex function $f:\R^n\to\R$ whose subdifferential satisfies \eql{int}{F(x)\subset\partial f(x)\te{for all}x\in\R^n.}

If $f:\R^n\to\R$ is convex, the subdifferential $\partial f:\R^n\mto\R^n$ satisfies $f(z+d)-f(z)\geq \partial f(z)^\T d$ for all $z,d\in\R^n$. If \r{int} is valid,
we hence conclude
(after the change of variables $z=x+u$ and $d=-u$) that $F$ is linked to $f$ by
\eql{disf}{
f(x+u)-f(x)\leq F(x+u)^\T u\te{for all}x,u\in\R^n.
}
This just is the local dissipation inequality for the system
\eql{sysf}{
\si x=x+u,\ \ y\in F(x+u)
}
and the supply rate $S(u,y)=u^\T y$, with $f$ being a storage function. Hence \r{sysf} is passive.
Conversely, passivity of \r{sysf} implies \r{disf} for some $f:\R^n\to\R$. It is then elementary to see that $f$ is convex, and \r{int} holds by the mere definition of the subdifferential. This proves the following result.

\theorem{\label{Tpas}The mapping $F:\R^n\mto\R^n$ has a convex primitive iff the system $\si x=x+u$, $y\in F(x+u)$ is passive.}

Since the linear system $\si x=x+u$ is controllable, passivity is equivalent to cyclo-passivity by Theorem~\ref{Tdis}.
Now recall that $F$ is said to be cyclically monotone if
\eql{cym}{
\sum_{j=0}^m F(v_j)^\T(v_j-v_{j+1})\geq 0
}
holds for any choice of points $v_0,\ldots,v_m\in\R^n$ and with $v_{m+1}:=v_0$.
It is not hard to establish the following link (as shown in Section~\ref{SLc}).

\lem{\label{Lcyc} The system $\si x=x+u$, $y\in F(x+u)$
is cyclo-passive iff $F$ is cyclically monotone.}


This leads to Rockafellar's celebrated result that $F$ has a convex primitive iff it is cyclically monotone
(\cite{Roc66}), with a proof based on dissipativity theory.

Instead of exploring these links any further, we rather  exploit Theorem~\ref{Tpas} for the generation of integral quadratic constraints in robustness analysis.
As a first step, we note that \r{sysf} is dissipative w.r.t.\ $(u,y)\mapsto u^\T y$ iff this holds for the
system and the supply rate obtained after the pre-compensation  $u=-x+v$. This transformation leads to the system
$\si x=v$, $y\in F(v)$ with supply rate $(v,y)\mapsto y^\T (v-x)$. Thus,
passivity of \r{sysf} and dissipativity of
\eql{sys1}{
\si x=u,\ \ y=-x+u,\ \ S(u,y)=F(u)^\T y
}
are equivalent. Note that the convex integrability of $F$ can as well be expressed by dissipativity of \r{sys1}.

As a consequence, for some convex function $f:\R^n\to\R$,
\eql{sysp}{
\si x=u\ \ y=-x+u,\ \ S(u,y)=\partial f(u)^\T y
}
is dissipative as certified by the storage function $f$. The linear system or filter in \r{sysp} operates on $u$, and the supply rate is the inner product of the filter's output and $\partial f(u)$.
Interestingly, dissipativity is as well guaranteed if the filter operates on $\partial f(u)$ and the supply rate
is the inner product of $y$ and $u$, which means to consider
\eql{sysd}{\si x=v,\ \ y=-x+v,\ \ v\in\partial f(u),\ \ S(u,y)=y^\T u.}
A suitable storage function is the Fenchel conjugate of $f$, which is defined as
$$
f^*(x)=\sup_{w\in\R^n}[x^\T w-f(w)]\te{for}x\in\R^n.
$$
Note that $f^*$ takes finite values on the image of the subdifferential mapping $\partial f$ given by
$$
X_{\partial f}:=\{x\in\R^n \ |\  x\in \partial f(u)\te{for some}u\in \R^n\}.
$$
Our findings are summarized as follows and proved in Section~\ref{STdis2}.
\theorem{\label{Tdis2}Let $f:\R^n\to\R$ be convex and $U=Y=\R^n$. Then \r{sysp} is dissipative with storage function $f$
on $X=\R^n$.
Moreover,  \r{sysd} is dissipative with storage function $f^*$ on the state-set $X^*:=X_{\partial f}$.
}

These results extend to longer shifts, if we replace $\si$ by $\si^k$ for $k=1,\ldots,\nu\in\N$ in the dynamics of \r{sysp} and \r{sysd}. First order representation of these dynamics are
described with
\eql{jor}{
\Afk:=\Af\otimes I_n,\ \ \Bfk:=\Bf\otimes I_n,\ \ \Cfk:=\Cf\otimes I_n,
}
where $A_\nu\in \R^{\nu\times\nu}$ is the upper Jordan block with eigenvalue zero, $B_\nu:=e_\nu\in\R^\nu$, $C_{k,\nu}:=-e_k^\T\in\R^\nu$, $I_n$ is the identity matrix of dimension $n$, and $\ot$ is the Kronecker product.

\corollary{\label{Cdis}Let $f:\R^n\to\R$ be convex, $k\in\{1,\ldots,\nu\}$ and $U=Y=\R^n$. Then
\eql{syspk}{
\si x=\Afk x+\Bfk u,\ \ y=\Cfk x+u,\ \ S(u,y)=\partial f(u)^\T y
}
is dissipative on $X:=\times_{j=1}^\nu\R^{n}$ with the storage function $V_k(x)=\sum_{j=k}^\nu f(x^j)$ for
$x=(x^1,\ldots,x^\nu)\in X$.
Moreover,
\eql{sysdk}{
\si x\in\Afk x+\Bfk v,\ \ y=\Cfk x+v,\ \ v\in\partial f(u)
}
is passive on $X^*:=\times_{j=1}^\nu X_{\partial f}$
with the storage function
$V_k^*(x)=\sum_{j=k}^\nu f^*(x^j)$ for $x=(x^1,\ldots,x^\nu)\in X^*$.
}

The proof is given in Section~\ref{SCdis}.
Corollary~\ref{Cdis} forms the basis for generating causal and anti-causal OZF multipliers in Section~\ref{Szf}.
Let us emphasize that the construction of anti-causal multipliers involves Fenchel conjugation.

\section{Strict Dissipativity and Linear Systems}\l{Sdis1}

We now briefly recap strict dissipativity theory for linear systems.
A detailed discussion for continuous-time systems can be found in the recent survey paper by \cite{Sch22}.

\definition{
The nonlinear system \r{sys} is strictly dissipative with respect to $S:U\times Y\mto\R$ if there exists a storage function $V:X\to\R$
and some $\eps>0$ such that
\eql{di}{
V(x_{t_2})\leq V(x_{t_1})+\sum_{t=t_1}^{t_2-1} S(u_t,y_t)-\eps\sum_{t=t_1}^{t_2-1}\left(\|x_t\|^2+\|u_t\|^2\right)
}
holds for all admissible trajectories and all time instances $t_1,t_2\in\N_0$ with $t_1\leq t_2$.}

For a linear system
\eql{lsys}{\si x=Ax+Bu,\ \ y=Cx+Du}
and a homogeneous quadratic supply rate
\eql{qsu}{S_P(u,y):=\mat{c}{y\\u}^\T\!\! P\mat{c}{y\\u},\ \ P=P^\T \in\R^{(k+m)\times (k+m)},}
strict dissipativity with a general storage function is equivalent to strict dissipativity with a homogeneous quadratic storage function. This leads to the following key result.
\theorem{\label{Tsdi}
The linear system \r{lsys} is strictly dissipative with respect to the quadratic supply rate \r{qsu} iff the there exists a symmetric solution $X$ of the strict LMI
\eql{lmi}{\mat{cc}{A&B\\I&0}^\T \!\mat{cc}{X&0\\0&-X}\mat{cc}{A&B\\I&0}\cl
\mat{cc}{C&D\\0&I_m}^\T \!\!P\mat{cc}{C&D\\0&I_m}.}
If $\eig(A)\cap\partial\Dset=\emptyset$ and with $G(\z):=C(\z I-A)^{-1}B+D$,
this is equivalent to the frequency-domain inequality
\eql{fdi}{
0\cl \mat{c}{G(\la)\\I_m}^*P\mat{c}{G(\la)\\I_m}\te{for all}\la\in\partial\Dset.
}
}
In summary, strict dissipativity can be assured by checking the feasibility of the LMI \r{lmi} or, under a mild assumption, through the frequency domain inequality \r{fdi}.

\section{Robust Stability Analysis}\l{Srob}

Let us now consider the feedback interconnection of
\eql{sy}{
x_{t+1}=\As x_t+\Bs w_t,\ z_t=\Cs x_t+\Ds w_t
}
with the subgradient nonlinearity
\eql{nl}{w_t\in\partial f(z_t)}
for $t\in\N_0$. We assume that $\As\in\R^{n\times n}$ is Schur and that $f:\R^d\to\R$ is convex with
$f(0)=0$ and $0\in\partial f(0)$.
Hence $f$, $f^*$ satisfy $f(x)\geq 0$ and $f^*(x)\geq f^*(0)=0$ for all $x\in\R^d$.

The current interest in such classical Lur'e system has been emphasized in the introduction, see also Section~\ref{Sdis}.

The interconnection \r{sy}-\r{nl} is said to be stable if
there exists some $L>0$ such that
\eql{sta}{
\sum_{t=0}^{\infty}\left(\|x_t\|^2+\|w_t\|^2\right)\leq L\|x_0\|^2
}
holds for all trajectories of \r{sy}-\r{nl}.
Now note, by \r{disf} for $F=\partial f$,  that the static system \r{nl}
is dissipative w.r.t. $S(z,w)=z^\T w$, i.e., it is passive.
Therefore, the classical passivity theorem guarantees
stability of \r{sy}-\r{nl} if the linear system \r{sy}
is strictly dissipative with respect to the supply rate $S(w,z)=-w^\T z$ (see, e.g.,
\cite{Bro04,BroTan20} in continuous-time).

This passivity-based stability test can be substantially improved by imposing dissipativity
constraints after passing the input and output signals of \r{sy} or \r{nl} through some dynamic filter
\eql{fil}{
\si\xi=A_\Psi\xi+B_{\Psi}\mat{c}{z\\w},\ v=C_\Psi\xi+D_{\Psi}\mat{c}{z\\w},\ \xi_0=0
}
with $\xi$ of dimension $n_\Psi$.
With the state $\eta=\col(\xi,x)$ of dimension $n_\Psi+n$, the interconnection of \r{fil} with \r{sy} is
\eql{filsy}{
\si\eta=\Afs\eta+\Bfs w,\ \ v=\Cfs\eta+\Dfs w,\te{where}
}
$\tiny\arraycolsep0.19ex\Afs\!=\!
\mat{cc}{
A_\Psi&B_{\Psi}\smat{c}{\Cs\\0}\\
0     &\As         },\
\Bfs\!=\!
\mat{cc}{
B_{\Psi}\smat{c}{\Ds\\I}\\B},\
\Cfs\!=\!\mat{cc}{C_\Psi&D_{\Psi}\smat{c}{\Cs\\0}},\ \Dfs\!=\!
D_{\Psi}\smat{c}{\Ds\\I}.
$

Let us now formulate a discrete-time robust stability result whose continuous-time counterpart has been first presented
in \cite{SchVee18}.

\theorem{\label{Trs}
Let $P=P^\T$ and suppose the following holds.
\enu{
\item 
All trajectories of \r{fil}
under the constraint \r{nl} satisfy the following IQC with a quadratic terminal cost defined by some matrix
$Z=Z^\T$:
\eql{iqcZ}{
\sum_{t=0}^{T-1} v_t^\T Pv_t-\xi_{T}^\T Z\xi_{T}\geq 0\te{for all}T\in \N.
}
\item The system \r{filsy} is strictly dissipative w.r.t. the supply rate $S(w,v)=-v^\T Pv$ as certified by $\c{X}=\c{X}^\T$.
\item The certificates $\c{X}=\tiny \arraycolsep.2ex\mat{cc}{X_{\Psi}&W\\W^\T &X}\in\R^{(n_\Psi+n)\times(n_\Psi+n)}$ and $Z\in\R^{n_\Psi\times n_\Psi}$ are coupled as
\eql{pos}{
\mat{cc}{X_{\Psi}+Z&W\\W^\T &X}\cg 0.
}
}

Then there exists some $\eps>0$ such that all the trajectories of \r{sy}-\r{nl} satisfy
\eql{sta}{
x_T^\T Yx_T+\eps\sum_{t=0}^{T-1}\left(\|x_t\|^2+\|w_t\|^2\right)\leq x_0^\T X x_0
}
for all $T\in\N$ where $Y:=X-W^\T(X_{\Psi}+Z)^{-1}W$. In particular,
the loop \r{sy}-\r{nl} is stable.
}

The proof is given in Section~\ref{STrs}.
Let us contrast Theorem~\ref{Trs} with a stability result
based on soft IQCs due to \cite{MegRan97},
which can be proved in discrete-time as in \cite{Sei15} and \cite{VeeSch14} without relying on homotopy arguments.
This requires that \r{sy}-\r{nl} is well-posed in the sense defined in these references. Moreover,
$f$ needs to be differentiable with $\|\nabla f(x)\|\leq M\|x\|$ for all $x\in\R^n$ and some $M\geq0$.
Finally, $A_\Psi$ of the filter \r{fil}
with the transfer matrix $\Psi(\z)=C_\Psi(\z I-A_\Psi)^{-1}B_\Psi+D_\Psi$
is assumed to be stable.

\theorem{\label{Tiqc}
Let $P=P^\T$ and suppose the following holds.
\enu{
\item All trajectories of \r{fil} with $z\in l_2^d$ and $w=\nabla f(z)$  satisfy the soft IQC
$\sum_{t=0}^\infty v_t^\T Pv_t\geq 0.$
\item The system \r{filsy} is strictly dissipative with respect to the supply rate $S(w,v)=-v^\T Pv$.
\item The left-upper/right-lower $d\times d$-block of
$\Psi(\la)^*P\Psi(\la)$ is positive/negative semi-definite for all $\la\in\partial\Dset$.
}
Then the loop \r{sy}-\r{nl} is stable.
}

In both results, it is a misnomer to talk about an integral quadratic constraint in a), since the terminology results from the continuous-time counterparts involving integration instead of summation. Moreover,
the transfer matrix $\Psi(1/\z)^TP\Psi(\z)$ is referred to as a dynamic multiplier if
$\Psi(\z)$ is not constant, and it is said to be static otherwise.


Observe that the IQC in Theorem~\ref{Tiqc} a) is called soft since it is formulated
on the infinite time-horizon. Although the conditions in b) are equivalent, Theorem~\ref{Tiqc} does not involve any sign-constraint on the corresponding certificate.
For these reasons, Theorem~\ref{Tiqc}
does not allow to draw conclusions about the transient behaviour of the system state.
In contrast, Theorem~\ref{Trs} leads to guaranteed point-wise in time constraints for the future state-trajectory  as in \r{sta}. In classical dissipation-based stability tests, this is achieved by taking $Z=0$ in Theorem~\ref{Trs} a), which leads to a so-called hard IQC. Then Theorem~\ref{Trs} c) requires that
the certificate $\c{X}$ is positive definite, which involves possibly severe conservatism
(see Section~\ref{Salg}).

A much more detailed discussion of these links and further consequences  are provided in \cite{SchVee18} and \cite{Sch22}.

\section{O'Shea-Zames-Falb Multipliers}\label{Szf}

We now establish how to construct IQCs with a nontrivial terminal cost $Z\neq 0$ as appearing in  Theorem~\ref{Trs} a).
To this end, let $y=\psi_{k,\nu}(u)$ be defined by the filter
$$
\mat{c}{\si x\\y}=\mat{cc}{\Af&\Bf\\\Cf&\Df}\mat{c}{x\\u}\te{with}x_0=0.
$$
\begin{figure}[]
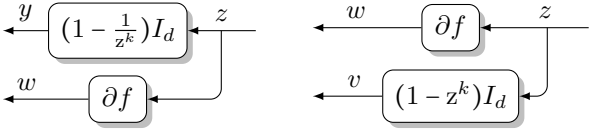
\center
\fbzf\hspace{.5cm}\fbzfa
\caption{Basic causal and non-causal OZF multipliers.\label{fig0}}
\end{figure}

Since $f$ is nonnegative and if
$k\in\{1,\ldots,\nu\}$, Corollary~\ref{Cdis} implies for $z\in l_{2e}^d$ and
$y=(\psi_{k,\nu}\ot I_d)(z)$ with $w\in\partial f(z)$ that the hard IQC
$\sum_{t=0}^{T-1} w_t^\T y_t\geq 0$ is satisfied for all $T\in\N$. Since
$\psi_{k,\nu}$ has the causal transfer function $1-1/\z^k$,
this IQC is formulated for the configuration on the left in Fig.~\ref{fig0}
involving a causal filter.
Note that $\psi_{0,\nu}(\z)=1$.  Since $\partial f$ is passive, we also get a
static hard IQC for $k=0$.

By Corollary~\ref{Cdis},
the trajectories of $\t y=(\psi_{k,\nu}\ot I_d)(w)$ with $w\in\partial f(z)$ satisfy
$\sum_{t=0}^{T-1}z_t^\T \t y_t \geq 0$ for all $T\in\N$ and $k\in\{0,1,\ldots,\nu\}$.
This can be interpreted as an IQC for the configuration on the right in
Fig.~\ref{fig0}. This involves a non-causal filter, thus defining a non-causal OZF-multiplier.

Let us now conically combine all these IQCs. To this end,
fix possibly different values of $\nup,\nud\in\N_0$ and define
\eql{filpd}{
\psi_\nu(\la):=\sum_{k=0}^\nu \lap_k\psi_{k,\nu}\te{for}\la:=(\la_0,\ldots,\la_\nu)\in\R^{\nu+1}.
}
Note that $\psi_\nu$ admits a state-space realization in terms of
$(\Afp,\Bfp,\Cfp,\Dfp)$ with
$C_\nu(\lap):=\sum_{k=1}^\nu \lap_k C_{k,\nu}$ and $D_\nu(\lap):=\sum_{k=0}^\nu\lap_\nu.$
If the vectors $\la\in\R^{\nup+1}$ and $\t\la\in\R^{\nud+1}$ are nonnegative, we then obtain
for any $z\in l_{2e}^d$, $w\in \partial f(z)$ and
$y=(\psi_\nu(\la)\otimes I_d)(z)$, $\t y=(\psi_{\t\nu}(\t\la)\ot I_d)(w)$ the hard IQC
\eql{zf3}{\sum_{t=0}^{T-1}\left(w_t^\T y_t+z_t^\T \t y_t\right)\geq 0\te{for all}T\in\N.}
This is condition a) in Theorem~\ref{Trs} for $Z=0$. Indeed, with
\eql{f1}{
\arr{rclrcl}{
A_\Psi&:=&\mat{cc}{\Afp&0\\0&\Afd}\ot I_d,\ \ & B_\Psi&:=&\mat{cc}{\Bfp&0\\0&\Bfd}\ot I_d,\\
C_\Psi&:=&\mat{cc}{I_{\nup+\nud}\\0}\ot I_d,\ \ &
D_\Psi&:=&\mat{cc}{0\\I_{2}}\ot I_d,}
}
and the symmetric matrix 
\eql{f3}{
P(\lap,\lad):=
\mat{cc|cc}{0&0&0&\bullet\\0&0&\Cfd^\T&0\hl&&\\[-2ex] 0&\Cfd  &0&\bullet\\ \Cfp &0&\Dfp+\Dfd^\T&0}\ot I_d,
}
the response of \r{fil} to $z\in l_{2e}^d$ and $w\in\partial f(z)$ satisfies
\eql{h0}{v_t^\T P(\lap,\lad)v_t=(w_t^\T y_t+z_t^\T \t y_t)+(y_t^\T w_t+\t y_t^\T z_t)}
for all $t\in\N_0$. By \r{zf3}, we infer a) in Theorem~\ref{Trs} for $Z=0$.
To generate IQCs with $Z\neq 0$, let us introduce the lifted state-space representations of the filters
$y=\psi_\nup(\lap)(z)$ and $\t y=\psi_\nud(\lad)(w)$ as in the notation section. Since the filter's initial conditions are zero, this results in
\eql{fill}{
\mat{c}{x_T\\y^T}=\mat{c}{\Bfp^T\\\Dfpl{T}}z^T,\ \mat{c}{\t x_T\\\t y^T}=\mat{c}{\Bfd^T\\\Dfdl{T}}w^T.
}
We will show that \r{iqcZ} indeed holds with the filter \r{f1} and the matrix $P(\lap,\lad)$ in \r{f3} for the terminal cost matrix
\eql{Z}{
Z(E):=\mat{cc}{0&E^\T\ot I_d\\E\ot I_d&0},
}
if the coefficient vectors $\lap$, $\lad$ and the free matrix $E$ render
\eql{MT}{
M_T(\lap,\lad,E):=\Dfpl{T}+\Dfdl{T}^\T-(\Bfd^{T})^\T E\Bfp^{T}
}
doubly hyperdominant for $T=\nup+\nud+1$.



\theorem{\label{Tzf}Suppose that $\lap\in\R^{\nup+1}$, $\lad\in\R^{\nud+1}$ and $E\in\R^{\nud\times\nup}$ are chosen with
$$
M_{\nup+\nud+1}(\lap,\lad,E)\in\c{H}^{(\nup+\nud+1)\times(\nup+\nud+1)}.
$$
Then the nonlinearity $\partial f$ satisfies the IQC \r{iqcZ} with the filter defined by \r{f1} for the supply rate matrix $P(\lap,\lad)$ in \r{f3} and the terminal cost matrix $Z(E)$ in \r{Z}.}

This is the second main result of this paper, whose proof is found in Section~\ref{STzf}. No analogous result is available in continuous-time.

%

\section{A Concrete Algorithm}\l{Salg}

Theorems~\ref{Trs} and \ref{Tzf} form the basis for generating various robust performance results with dynamic IQCs by known dissipativity arguments. A collection of such results can be found, e.g., in \cite{FetSch17d} and \cite{Sch22}.


For example, if $e=C_ex$ is an output of the interconnection \r{sy}-\r{nl}, we can target at the computation of some
tight $\ga>0$ such that the amplitude bound $\sup_{t\in\N}\|e_t\|\leq\ga$ holds for all initial conditions $x_0$ in the unit ball of $\R^n$ and all trajectories of \r{sy}-\r{nl}.

\corollary{\label{Cal}
The loop \r{sy}-\r{nl} is stable and all its trajectories satisfy
$$\sup_{t\in\N}\|C_ex_t\|\leq\ga_*(\nup,\nud)\|x_0\|$$
if $\ga_*(\nup,\nud)$ is determined as follows:
\enu{
\item Choose the dimensions $\nup,\nud\in\N_0$ to fix the multiplier complexity and set $T_0:=\nup+\nud+1$.
\item With the filter matrices \r{f1} construct \r{fil} and \r{filsy}.
\item With $P(.)$, $Z(.)$, $M_{T_0}(.)$ defined by \r{f3}, \r{Z}, \r{MT} and
the variables $\c{X}=\c{X}^\T$, $H=H^\T$, $E\in\R^{\nud\times\nup}$, $\lap\in\R^{\nup+1}$, $\lad\in\R^{\nud+1}$, $\ga\in\R$, introduce the LMIs
\eql{lmip}{
\arr{c}{
\mat{cc}{\Afs&\Bfs\\I&0}^\T \!\mat{cc}{\c{X}&0\\0&-\c{X}}\mat{cc}{\Afs&\Bfs\\I&0}+\hspace{2cm}\\\hspace{2cm}+
\mat{cc}{\Cfs&\Dfs}^\T \!\!P(\lap,\lad)\mat{cc}{\Cfs&\Dfs}\cl 0,\\
\mat{ccc}{H&0&C_e\\0&X_{\Psi}+Z(E)&W\\C_e^\T&W^\T &X}\cg 0,\ H\cl \ga I,\ X\cl \ga I,\\
M_{T_0}(\lap,\lad,E)\in\c{H}^{T_0\times T_0},}
}
where $\c{X}$ is partitioned as in Theorem~\ref{Trs}.
\item Then $\ga_*(\nup,\nud)$ denotes the infimum of all $\ga>0$ for which the LMIs \r{lmip} are feasible.
}
}
It is instrumental to observe that the constraints \r{lmip} are indeed affine in all variables. Hence, the computation of
$\ga_*(\nup,\nud)$ involves solving a standard semi-definite program.

\begin{figure}\center
\includegraphics[width=8.5cm]{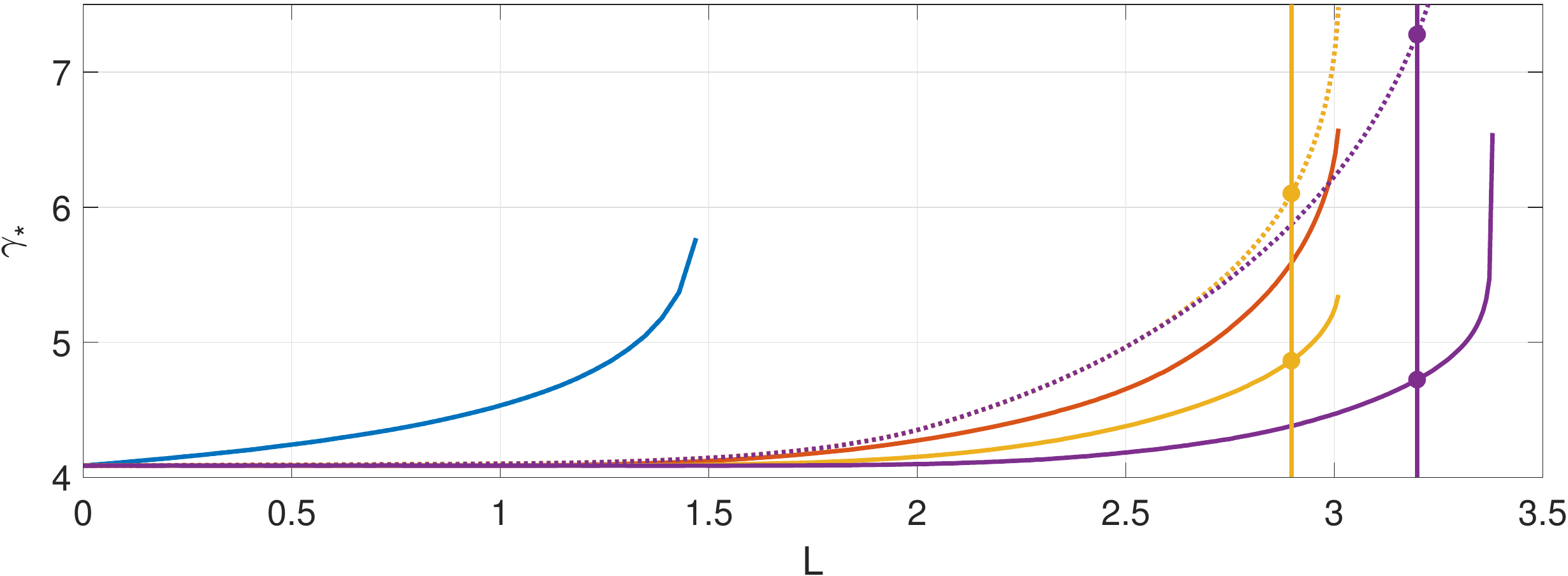}
\caption{Computed values of $\ga_*(0,0)$ (blue), $\ga_*(1,1)$ (red), $\ga_*(2,2)$ (yellow), and $\ga_*(3,3)$ (purple):
Full (dotted) lines are obtained with Corollary~\ref{Cal} (setting $E=0$).}\label{fig1}
\end{figure}

For a numerical example, take the system matrices
$$
\mat{c|c}{A&B\hl C&D_L\hl C_e&0}=\left(\begin{array}{ccccc|c}
0.1&1&0&0&0&0\\-0.24&0.1&-0.54&-0.35&0.84&0\\0&0&0.54&-0.24&0.59&0\\0&0&0&0.54&1&0\\0&0&0&-0.56&0.54&1.04
\hl -0.08&-0.17&0.13&0.09&-0.21&-1/L\hl
2&1.3&-1&-1.3&1.3&0
\end{array}\right)
$$
depending on $L\in(0,3.5]$, and compute $\ga_*(\nup,\nud)$ as in Corollary~\ref{Cal} for
$\nup=\nud\in\{0,1,2,3\}$ using the LMI-solver of \cite{ML20b} and Yalmip (\cite{Lof04}. The results are plotted over $L$ in Fig.~\ref{fig1} as full lines. They are compared with the values for hard IQCs (dotted curves), which are obtained with $E=0$ in Corollary~\ref{Cal}. The blue curve indicates the severe conservatism of results based on static IQCs ($\nup=\nud=0$).  For dynamic IQCs ($\nup=\nud>0$), the plots permit to quantify the conservatism
of hard IQCs if compared to those with a nontrivial terminal cost as newly developed in this paper. For example, the yellow and purple dots in Fig.~\ref{fig1} indicate a reduction of conservatism by about $20\%$ and $35\%$, respectively.

\section{Discussion and a Conjecture}\label{Sdis}

Due to the modularity of results based on dissipativity, we emphasize that Theorems~\ref{Trs} and \ref{Tzf}
can be seamlessly merged with the approaches in, e.g., \cite{YinSei20a} and \cite{PauGra21} to reduce
the conservatism in the determination of stability margins of neural network controllers.
Our simple numerical example demonstrates the striking potential limitation of static and hard dynamic IQCs, as employed in a whole stream of recent papers revolving around the generation of safety guarantees for static maps or dynamic systems involving neural networks.

In an another direction, the analysis of gradient descent algorithms for $m$-strongly convex and $L$-smooth functions
is addressed by \cite{LesRec16} based on causal OZF multipliers, while the extension to general OZF multipliers is proposed by \cite{MicSch20}. In the latter situation,
our results permit to guarantee transient properties of algorithms, next to the typically investigated
global exponential
stability. 

We conclude the paper with the conjecture that robustness analysis based on Theorems~\ref{Trs} and \ref{Tzf}
does not involve any more conservatism than robustness analysis based on Theorem~\ref{Tiqc}
with general OZF multipliers.


\section{Conclusions}\label{S5}

We have developed a novel discrete-time absolute stability result by dissipation techniques and
based on the notion of integral quadratic constraints with a nontrivial terminal cost. The benefit over
existing results has been shown by a numerical example, and we conjecture that it
is a lossless extension of a known robust stability result relying on O'Shea-Zames-Falb multipliers and soft IQCs.
The suggested link of dissipativity theory with convex analysis might offer avenues for systematically developing new
classes of multipliers in IQC based stability results.

{

}

\appendix
\section{Proofs}
\subsection{Proof of Lemma~\ref{Lcyc}} \label{SLc}

Let \r{sysf} be cyclo-passive. If we are given $v_0,\ldots,v_m,v_{m+1}=v_0\in\R^n$ and any $w_j\in F(v_j)$, we can take
$$
x_0=v_0\te{and}u_t=v_{m-t}-v_{m+1-t}\te{for}t=0,\ldots,m
$$
to generate a trajectory of \r{sysf} with
$$
x_t:=v_{m+1-t}\te{and}x_t+u_t=v_{m-t} \te{for}t=0,\ldots,m+1.
$$
With $y_t:=w_{m-t}$ we get
$y_t^\T u_t=w_{m-t}^\T (v_{m-t}-v_{m-t+1})$ for $t=0,\ldots,m+1$, and thus
\equ{\label{cyc}
 \sum_{j=0}^m w_j^\T(v_j-v_{j+1})=\sum_{t=0}^m y_t^\T u_t.
}
If noting $y_t\in F(v_{m-t})=F(x_t+u_t)$ and $x_0=v_{m+1}=v_0=x_{m+1}$, cyclo-passivity of \r{sysf} allows us to conclude
$\sum_{t=0}^m y_t^\T u_t\geq 0$, which implies \r{cym} due to \r{cyc}.

Conversely, suppose that $F$ is cyclically monotone. Let us pick any trajectory of \r{sysf}
with $x_0=x_{m+1}$. If defining
$$
v_j:=x_{m+1-j}\te{and}w_j:=y_{m-j}\te{for}j=0,\ldots,m+1,
$$
we have $w_j^\T(v_{j}-v_{j+1})=y_{m-j}^\T(x_{m-j+1}-x_{m-j})=y_{m-j}^\T u_{m-j}$
for $j=0,\ldots,m+1$, giving again \r{cyc}.
Since $w_j=y_{m-j}\in F(x_{m-j}+u_{m-j})=F(x_{m-j+1})=F(v_j)$ and
due to \r{cym}, we infer that $\sum_{t=0}^m y_t^\T u_t\geq 0$. Since the choice
of the round trip trajectory of \r{sysf} was arbitrary,  \r{sysf} is cyclo-passive.

\subsection{Proof of Theorem~\ref{Tdis2}}\label{STdis2}

We only need to prove the second statement.
Pick any $u\in U$ and $v\in\partial f(u)$. Then $0$ is a subgradient of $w\mapsto v^\T w-f(w)$ at $w=u$. By Fermat's principle, $w=u$ is optimal for the maximization in the definition of  $f^*(v)$, which implies
\eql{Ah1}{f^*(v)=v^\T u-f(u)<\infty.}
To show the dissipation inequality, take any $x \in X^*$ and pick $\t u\in\R^n$ with $x\in\partial f(\t u)$.
Then the value of $f^*(x)$ is finite, and
the definition of $f^*(x)$ implies $-f^*(x)\leq f(u)-x^\T u$.
Addition to \r{Ah1} gives $f^*(v)-f^*(x)\leq u^\T (v-x).$ We have proven that
$$f^*(v)-f^*(x)\leq u^\T (v-x)\te{for all}x\in X^*,\ u\in U,\ v\in\partial f(u).$$
This is the local dissipation inequality
for \r{sysd} on $X^*\times U$.\\
\mbox{}\hfill 

\subsection{Proof of Corollary~\ref{Cdis}}\label{SCdis}

If $(x,u)\in X\times U$, then  $\Afk x+\Bfk u=\col(x^2,\ldots,x^{\nu},u)$ and  $\Cfk x+u=u-x^k$ just by \r{jor}. This shows
\mun{V_k(\Afk x+\Bfk u)-V_k(x)=\sum_{j=k+1}^{\nu} f(x^j)+f(u)-\sum_{j=k}^\nu f(x^j)=\\=
f(u)-f(x^k)\leq \partial f(u)^\T (u-x^k)=\partial f(u)^\T (\Cfk x+u),}
where the inequality follows from the dissipation inequality for \r{sysp} with storage function $f$. This proves the first claim.

Similarly, $(x,u)\in X^*\times U$ implies that there exist $u_j\in\R^n$ with
$x\in\col(\partial f(u_1),\partial f(u_2),\ldots,\partial f(u_\nu))$. For any $u\in U$ and $v\in\partial f(u)$, we get
$\Afk x+\Bfk v\in\col(x^2,\ldots,x^\nu,v)$, thus
\mun{V_k^*(\Afk x+\Bfk v)-V_k^*(x)=
\!\!\!\sum_{j=k+1}^{\nu} \!\!\!f^*(x^j)+f^*(v)-\sum_{j=k}^\nu f^*(x^j)=\\=
f^*(v)-f^*(x^k)\leq u^\T (v-x^k)=u^\T (\Cfk x+v),}
by relying on the dissipation inequality for \r{sysd} with the storage function $f^*$. This completes the proof.

\subsection{Proof of Theorem~\ref{Trs}}\label{STrs}
Pick any trajectory of the loop \r{sy}-\r{nl}. By the dissipativity condition b) and Theorem~\ref{Tsdi}, there exists some $\eps>0$ such that the signals $z$ and $w$ of the loop trajectory filtered by \r{fil} fulfill the dissipation inequality
\mun{
\mat{c}{\xi_{T}\\x_{T}}^\T \!\!\mat{cc}{X_\Psi&W\\W^\T &X}\mat{c}{\xi_{T}\\x_{T}}+\sum_{t=0}^{T-1} v_t^\T Pv_t+\\+\eps\sum_{t=0}^{T-1} \left(\|\xi_t\|^2+\|x_t\|^2+\|w_t\|^2\right)\leq
\mat{c}{0\\x_0}^\T \!\!\mat{cc}{X_\Psi&W\\W^\T &X}\mat{c}{0\\x_0}
}
for all $T\in\N_0$. By \r{nl}, we can make use of \r{iqcZ} to infer
\mun{
\mat{c}{\xi_{T}\\x_{T}}^\T \!\!
\mat{cc}{X_\Psi+Z&W\\W^\T &X}\mat{c}{\xi_{T}\\x_{T}}
+\eps\sum_{t=0}^{T-1} \left\|\mat{c}{x_t\\w_t}\right\|^2\leq x_0^\T Xx_0
}
for all $T\in\N_0$. With c) and the standard fact
$$
x_T^\T Y x_T=\inf_{\xi\in\R^{n_\Psi}} \mat{c}{\xi\\x_{T}}^\T \!\! \mat{cc}{X_\Psi+Z&W\\W^\T &X}\mat{c}{\xi\\x_{T}},
$$
the proof is completed.
\mbox{}\hfill

\subsection{Proof of Theorem~\ref{Tzf}}\label{STzf}

We drop the arguments of $P(.)$, $Z(.)$ and $M_T(.)$ in \r{f3}, \r{Z} and \r{MT}, respectively, to lighten the notation. Let us first show
$M_T\in\c{H}^{T\times T}$ for $T\in\N$ different from $T_0:=\nup+\nud+1$.

Indeed, the cancellation of the first column of $\Bfp^{T_0}/\Bfd^{T_0}$ leads to $\Bfp^{T_0-1}/\Bfd^{T_0-1}$, while
canceling the first row and column of $\Dfpl{T_0}/\Dfdl{T_0}$ generate $\Dfpl{T_0-1}/\Dfdl{T_0-1}$, respectively.
Since principal sub-matrices of doubly hyperdominant matrices are double hyperdominant and
since $M_{T_0}$ is doubly hyperdominant by assumption, we conclude
that $M_{T_0-1}$ is doubly hyperdominant as well. This argument can be repeated to prove
the claim for $1\leq T<T_0$.

Let us now show the claim for $T=T_0+1$, while the general case follows by induction. By exploiting the particular structure of
$(\Afp,\Bfp,\Cfp,\Dfp)$, it is elementary to verify that $M_{T_0+1}$ looks with
$s=\sum_{k=0}^\nup\lap+\sum_{k=0}^\nud\lad$ as
\mun{
\mat{c:cccc|ccccccc}{
s           &-\lad_1       &\cdots    &-\lad_\nud &0          &0          &0         &\cdots &0          \hdl
-\lap_1     &s             &-\lad_1   &\cdots     &-\lad_\nud &0          &0         &\cdots &0          \\
-\la_2      &-\lap_1       &s         &-\lad_1    &\cdots     &-\lad_\nud &0         &\cdots &0          \\
\vdots      &\vdots        &\ddots    &\ddots     &\ddots     &           &\ddots    &\ddots &\vdots     \\
-\la_\nu    &\cdots        &-\lap_2   &-\lap_1    &s          &-\lad_1    &\cdots    &-\lad_\nud &0          \\
0           &-\lap_\nu     &\cdots    &-\lap_2    &-\lap_1    &s          &-\lad_1   &\cdots &-\lad_\nud \hl
\vdots      &\ddots        &\ddots    &\ddots     &\ddots     &\ddots     &\ddots    &\ddots &\vdots     \\
0           &\cdots        &0         &-\lap_\nu  &\cdots     &-\lap_2    &-\lap_1   &s      &-\lad_1    \\
0           & 0            &\cdots    &0          &-\lap_\nu  &\cdots     &-\lap_2   &-\lap_1&s          }-\\-
\mat{ccc}{
0&\cdots&0\hdl
0&\cdots&0\\
0&\cdots&0\\
\vdots&&\vdots\\
0&\cdots&0\\
0&\cdots&0\hl
1&\cdots&0\\
\vdots&\ddots&\vdots\\
0&\cdots&1
}E
\mat{c:cccc|ccccccc}{
0           & 0        &\cdots    &0         &0   &1       &0   &\cdots   &0\\
0           & 0        &\cdots    &0         &0   &0       &1   &\cdots   &0\\
\vdots      &\vdots    &          &\vdots    &\vdots         &\vdots  &\ddots  &\ddots   &\vdots\\
0           & 0        &\cdots    &0         &0   &0       &0   &\cdots   &1}
}
of dimension $(1+(\nu+1)+\nud)\times(1+(\nud+1)+\nu)$.
The dashed lines indicate where to find $M_{T_0}$ as a right-lower block, which is doubly hyperdominant by assumption.
As a consequence, $M_{T_0+1}$ has nonpositive off-diagonal elements. Moreover, the sum of the elements in the $(\nu+2)$-nd row ($(\nud+2)$-nd column) on top of the full horizontal line (left to the full vertical  line) is nonnegative.
All rows (columns) on top (to the left) of these are obtained by left (up) shifting this row (column)  and filling up the
last entries with zeros. This shows that the entries in the $i$-th row for $i=1,\ldots,\nu+1$ ($j$-th column for $j=1,\ldots,\nud+1$) sum up to nonnegative values as well. For these structural reasons, $M_{T_0+1}$ is double hyperdominant.

Now fix $T\in\N$. Due to \r{fill},  the lifted representations of
$y=(\psi_\nup(\lap)\ot I_d)(z)$ and $\t y=(\t\psi_\nud(\lad)\ot I_d)(w)$ read as
$$
\mat{c}{x_T\\y^T}=\mat{c}{\Bfp^T\ot I_d\\\Dfpl{T}\ot I_d}z^T,\ \mat{c}{\t x_T\\\t y^T}=\mat{c}{\Bfd^T\ot I_d\\\Dfdl{T}\ot I_d}w^T.
$$
If $z\in l_{2e}^d$ and $w\in\partial f(z)$, we get for these responses that
\mul{
\sum_{t=0}^{T-1}\left(w_t^\T y_t+z_t^\T \t y_t\right)-\t x_T^\T(E\ot I_d)x_T=\\=
(w^T)^\T (\Dfpl{T}\ot I_d)z^T+(z^T)^\T (\Dfdl{T}\ot I_d) w^T-\\\hspace{2cm}-(w^T)^\T (\Bfd^T\ot I_d)^\T (E\ot I_d)(\Bfp^T\ot I_d)z^T=\\=
(w^T)^\T \left( M_T\ot I_d\right)z^T.
\label{h1}
}
Note that this relation motivates the definition of $M_T$ after all.
Since $M_T$ is double hyperdominant and since
$$
w^T=\col(w_0,\ldots,w_{T-1})\in\col(\partial f(z_0),\ldots,\partial f(z_{T-1})),
$$
Lemma 2 in \cite{ManSaf05} guarantees that \r{h1} is nonnegative.

Let us finally recall \r{h0} for the response of the filter \r{fil} to the signals $z$ and $w$ and note that
$\xi_T=\col(x_T,\t x_T)$. For the matrix \r{Z} and together with \r{h1} added to its transposed version, we get
\mun{
\sum_{t=0}^{T-1}v_t^\T P(\nup,\nud)v_t-\xi_T^\T Z(E)\xi_T=\\=
\sum_{t=0}^{T-1}v_t^\T P(\nup,\nud)v_t-\t x_T^\T(E\ot I_d)x_T-x_T^\T(E^\T\ot I_d)\t x_T=\\=
(w^T)^\T \left( M_T\ot I_d\right)z^T+(z^T)^\T \left( M_T^\T\ot I_d\right)w^T\geq 0.
}
Since $T\in\N$ was arbitrary, the proof is concluded.

\subsection{Proof of Corollary~\ref{Cal}}\label{SCal}

Take any $\ga>\ga_*(\nup,\nud)\geq 0$. For this value of $\ga$, the LMIs \r{lmip} are feasible. Taking Schur complements gives
\eql{h2}{
\frac{1}{\ga}C_e^\T C_e\cle C_e^\T H^{-1}C_e\cl X-W^\T(X_{\Psi}+Z(E))^{-1}W
}
as a consequence of the second and third LMI. For any trajectory of the loop \r{sy}-\r{nl},
the fourth one implies $x_0^TX x_0\leq\ga\|x_0\|^2$.
On the one hand, Theorem~\ref{Trs} for $T\to\infty$ implies $\|x\|_2^2+\|w\|_2^2\leq \frac{\ga}{\eps}\|x_0\|^2$ for some $\eps>0$, which shows stability. On the other hand,
Theorem~\ref{Trs} combined with \r{h2} also guarantees
$$
\frac{1}{\ga}\|C_ex_T\|^2\leq \ga \|x_0\|^2\te{for all}T\in\N,
$$
which shows $\sup_{t\in\N}\|C_ex_t\|\leq\ga\|x_0\|$.
Since $\ga>\ga_*(\nup,\nud)$ was chosen arbitrarily, the proof is concluded.

\end{document}